\documentclass[a4paper,11pt,oneside  ]{article}	
\pdfoutput=1

\usepackage{a4wide,color,array}
\usepackage{cite}
\usepackage{hyperref}
\usepackage{amsmath,amssymb,amsthm,color}
\hypersetup{linktocpage,colorlinks, citecolor={magenta}}
\usepackage[english]{babel}
\usepackage[utf8]{inputenc}

\usepackage{xspace}
\usepackage{bm}				
\usepackage{bbm,upgreek}		
\usepackage[e]{esvect}		
\usepackage{esint}			
\usepackage{etoolbox}			
\usepackage{calc}				
\usepackage{eso-pic}			
\usepackage{datetime}			

\usepackage{lmodern}
\numberwithin{equation}{section}

\usepackage{enumitem}
\setlist{nosep}
\setlist{noitemsep}

\newcommand{\R}{\mathbf{R}}

\newtheorem{theorem}{Theorem}
\newtheorem{proposition}{Proposition}[section]
\newtheorem{lemma}[proposition]{Lemma}
\newtheorem{corollary}[proposition]{Corollary}
\newtheorem{remark}[proposition]{Remark}
\newtheorem{definition}[proposition]{Definition}

\theoremstyle{plain}

\theoremstyle{definition}

\newcommand{\cref}[1]{Corollary~\ref{c.#1}}

\def \1{\mathbf{1}} 

\def\({\left(}
\def\){\right)}

\def\XXint#1#2#3{{\setbox0=\hbox{$#1{#2#3}{\int}$}
		\vcenter{\hbox{$#2#3$}}\kern-.5\wd0}}

\renewcommand{\P}{\mathbb{P}}

\renewcommand{\tilde}{\widetilde}

\makeatletter
\def\namedlabel#1#2{\begingroup
	#2%
	\def\@currentlabel{#2}%
	\phantomsection\label{#1}\endgroup
}
\makeatother

\def \d {\mathsf{d}} 

\def \emp {\mathsf{emp}}  
\def \F {\mathsf{F}}  

\title{Emergence of a Poisson process in weakly interacting particle systems}
\author{David Padilla-Garza, Luke Peilen, Eric Thoma}
\date{}

\begin{document}
	
	\maketitle
	
	\maketitle
	
	\begin{abstract}
		We consider the Gibbs measure of a general interacting particle system for a certain class of ``weakly interacting" kernels. In particular, we show that the local point process converges to a Poisson point process as long as the inverse temperature $\beta$ satisfies $N^{-1} \ll \beta \ll N^{-\frac{1}{2}}$, where $N$ is the number of particles. This expands the temperature regime for which convergence to a Poisson point process has been proved.  
	\end{abstract}
	
	\section{Introduction and motivation}
	
	In this paper, we will be interested in a system of particles that interact via a pairwise interaction, and are confined by an external potential. We will consider a system of $N$ particles that lie in $d-$dimensional Euclidean space. This is modeled by the Hamiltonian 
	\begin{equation}
		\label{eq:hamiltonian}
		\mathcal{H}_{N}(X_{N}) = \sum_{i \neq j} g(x_{i} - x_{j}) + N \sum_{i=1}^{N} V(x_{i}),
	\end{equation}
	where $g: \mathbf{R}^{d} \to \mathbf{R}$ is the pair-wise interaction, $X_{N} := (x_{1}, ... x_{N}) \in (\mathbf{R}^d)^N$, and $V: \mathbf{R}^{d} \to \mathbf{R}$ is the confining potential. We will call this an interacting particle system. We will be interested in the behaviour of such a system for large but finite $N$.
	
	We are particularly interested in limiting behavior of the law of the local point process, which is defined as
	\begin{equation} \label{eq:empproc}
		\Xi := \sum_{i=1}^N \delta_{N^{1/d}(x_i - x^\ast)}
	\end{equation}
	for an artibrary centering point $x^\ast$. The scaling of $N^{1/d}$ will reflect that the particles $\{x_i\}_{i=1}^N$ typically occupy a volume of order $1$, so the local point process (in the weak topology) captures microscopic behavior of the system near $x^\ast$ as $N \to \infty$.
	
	Our main goal will be to prove that the law of $\Xi$ is asymptotically Poissonian in certain temperature scalings $\beta = \beta_N$ and for a broad class of ``weak interaction" kernels $g$. Before moving to precise statements in Section 2, we discuss more broadly the study of interacting particle systems with Hamiltonian of the form \eqref{eq:hamiltonian}.

	While the weakly interacting kernels that we consider require continuity, we are motivated by the study of Coulomb and Riesz interactions. A very frequent form of the interaction $g$ for $d \geq 2$ is given by the Coulomb kernel:
	\begin{equation}
		\label{Coulombcase}
		g(x)=
		\begin{cases}
			\frac{1}{|x|^{d-2}} &\text{ if } d \geq 3, \\
			- \log(|x|) &\text{ if } d = 2.
		\end{cases}
	\end{equation}
	Note that $g$ satisfies that 
	\begin{equation}
		\label{eq:fundamental}
		-\Delta g  = c_{d} \delta_{0},
	\end{equation}
	where $\delta_{x}$ denotes a Dirac delta at $x$, for all $d \geq 2$ and for some constant $c_{d}$ depending only on $d$.
	
	In $d=1$, a very frequent form of $g$ is given by 
	\begin{equation}
		\label{1dcase}
		g(x)=- \log(|x|) \text{ if } d = 1.
	\end{equation}
	This is the same formula as in $d=2$, but in dimension $1$ $g$ is not the fundamental solution of Laplacian (i.e. it does not satisfy equation \eqref{eq:fundamental}); instead it is a solution kernel for a fractional Laplacian.
	
	In dimension $d \geq 1$, a generalization of the Coulomb kernel is given by the Riesz kernel:
	\begin{equation}\label{Rieszcase}
		g(x) = \frac{1}{|x|^{d-2s}},
	\end{equation}
	with $0<s<\min\{ \frac{d}{2},1 \}$. In this case $g$ satisfies that
	\begin{equation}
		(- \Delta)^{s} g = c_{d,s} \delta_{0}.
	\end{equation}
	For some constant $c_{d,s}$ depending only on $d,s$. We will refer to this setting as the Riesz case.
	
	The study of minimizers of \eqref{eq:hamiltonian} when $g$ is either a Coulomb or Riesz interaction is an active field of research, with possible applications in approximation theory \cite{dragnev2016minimum,hardin2017generating,p2017next,reznikov2018covering,borodachov2018optimal,hardin2019local,dragnev2020constrained,chafai2022solution,chafai2023threshold}. However, in this paper we will take a different approach, and focus on the positive-temperature regime. In other words, we will look at the Gibbs measure associated to Hamiltonian \eqref{eq:hamiltonian}, instead of focusing only on minimizers. This Gibbs measure is given by:
	\begin{equation}
		\label{eq:Gibbs}
		\mathrm d \mathbb{P}_{N, \beta} (X_{N}) = \frac{1}{Z^{V}_{N, \beta}} \exp \left( - \beta \mathcal{H}_{N}(X_{N}) \right) \, \mathrm d X_{N},
	\end{equation}
	where
	\begin{equation}
		\label{eq:partfunc}
		Z^{V}_{N, \beta} = \int_{\mathbf{R}^{d \times N}} \exp \left( - \beta \mathcal{H}_{N}(X_{N}) \right) \, \mathrm d X_{N}
	\end{equation}
	is the partition function, and $\beta >0$ is the inverse temperature which may depend on $N$. \\
	
	The Gibbs measure \eqref{eq:Gibbs} is linked to Random Matrix Theory if $g$ is either a Coulomb, or a $1d$ log interaction. In $d=1$ and if $\beta=1,2$ or $4$ and if $V$ is quadratic, (\ref{eq:Gibbs}) gives the density of eigenvalues of matrices in the Gaussian Orthogonal, Unitary or Symplectic ensembles, respectively. These ensembles are well-studied because they are determinantal ($\beta=2$) or Pfaffian ($\beta=1,4$), which allows one to explicitly write down the correlation kernels. In the specific case $d=2, \beta=2, V(x)=|x|^{2}$, equation \eqref{eq:Gibbs} corresponds to the Ginibre ensemble, which also has a determinantal structure. Higher dimensional Coulomb gases, and Riesz gases are also an active field of research \cite{serfaty2023gaussian,armstrong2021local,leble2018fluctuations,leble2017large,petrache2017next,rougerie2016higher,boursier2021optimal,boursier2022decay, peilen2024log, thoma2024overcrowding}, and extending our techniques to the Coulomb gas is work in progress.

	At a macroscopic level, an interacting particle system is well-described by its empirical measure, defined as
	\begin{equation} \label{eq:emp}
		\mathrm{emp}_{N} := \frac{1}{N} \sum_{i=1}^{N} \delta_{x_{i}}. 
	\end{equation}
	The behaviour of the empirical measure is very different depending on the scaling of the temperature. If the temperature is not too large ($\frac{1}{N} 
	\ll \beta$), the empirical measure converges to the equilibrium measure, denoted $\mu_{V}$, and defined as the minimizer of the mean-field functional:
	\begin{equation} \label{eq:defmuV}
		\mu_{V} := \mathrm{argmin}_{\mu} \mathcal{E}(\mu) + \int_{\mathbf{R}^{d}} V \,  \mathrm d \mu,
	\end{equation}
	where
	\begin{equation}
		\mathcal{E}(\mu) := \iint_{\mathbf{R}^{d} \times \mathbf{R}^{d} } g(x-y) \,  \mathrm d \mu_{x} \,  \mathrm d \mu_{y}.
	\end{equation}
	This convergence happens in a large deviations sense \cite{garcia2019large}. This measure has compact support as long as $V$ grows fast enough at infinity (see section \ref{sect:prelims}). 
	
	On the other hand, for large temperatures ($\beta \simeq \frac{1}{N}$), the empirical measure does not converge to the equilibrium measure. Instead, the effect of temperature and entropy is large enough that particles are not confined to a compact set. In this case, the empirical measure converges to the thermal equilibrium measure, denoted $\mu_{\theta}$ and defined as 
	\begin{equation} \label{eq:defmutheta}
		\mu_{\theta} := \mathrm{argmin}_{\mu} \mathcal{E}(\mu) + \int_{\mathbf{R}^{d}} V \,  \mathrm d \mu + \frac{1}{\theta} \mathrm{ent}[\mu],
	\end{equation}
	where $\theta := N \beta$ and $\mathrm{ent}[\cdot]$ is the Boltzmann-Gibbs entropy:
	\begin{equation}
		{\rm ent}[\mu ]=
		\begin{cases}
			\int_{\mathbb{R}^{d}}  \log  ({ \frac{\mathrm d\mu}{\mathrm dx}})  \, \mathrm  d{\mu} \quad  \text{ if } \mu \text{ is absolutely continuous w.r.t. Lebesgue measure} \  \mathrm dx \\
			\infty\quad  \text{ o.w.}
		\end{cases}
	\end{equation} 
	Unlike the empirical measure, the thermal equilibrium measure is everywhere positive (see section \ref{sect:prelims}). In the case $g=0$ and $V(x)=|x|^{2}$, which corresponds to independent Gaussian particles, the thermal equilibrium measure is a Gaussian, with variance determined by $\theta$. Thus, the thermal equilibrium measure can be seen as an interpolation between a Gaussian probability density, and the equilibrium measure. 
	
	The empirical measure \eqref{eq:emp} and local point process \eqref{eq:empproc} exhibit very different behavior in the large $N$ limit. Due to equally weighting the $N$ particles, the empirical measure typically exhibits a law of large numbers effect and converges (weakly) to a deterministic limit. The local point process is sensitive to (and, in the weak topology, effectively determined by) each particle in a large microscopic neighborhood of the centering point. It exhibits random, temperature dependent behavior even in the $N \to \infty$ limit, except at very low temperature.
	
	At these low temperatures, the local point process is expected to be well approximated by minimizers of energies associated to the Gibbs measure Hamiltonian, an extremely delicate family of variational problems. For example, in certain dimensions, the local point process is conjectured to converge to a periodic lattice minimizing a next-order (renormalized) energy \cite{petrache2020crystallization,sandier20151d} in some generality. In dimensions $8$ and $24$, this conjecture has also been proved, in this case due to a link to the Cohn-Kumar conjecture \cite{petrache2020crystallization}. Given the difficulty in working with the local point process, it is often convenient to work instead with the empirical and tagged empirical fields, which are less refined but more manageable observables describing the microscopic behaviour of the particle system \cite{sandier20152d,rougerie2016higher,petrache2017next}. 
	
	On the other extreme, in the large temperature regime ($\beta \simeq \frac{1}{N}$), the law of the local point process completely thermalizes and is given by a Poisson process in the large $N$ limit. This was proved in the context of general integrable interactions in \cite{lambert2021poisson}, generalizing the work of \cite{benaych2015poisson,nakano2020poisson} for $\beta-$ensembles, which correspond to interacting particle systems with a Coulomb interaction. The temperature regime for Poissonian behavior is expected to be typically much wider, an expectation we confirm for our class of interactions.
	
	The microscopic behaviour of a particle system isn't necessarily asymptotically deterministic or Poissonian. For example, in the Coulomb case, if $\beta \simeq N^{1 - \frac{2}{d}}$ then the microscopic behaviour is stochastic, but governed by a probability measure that favours low-renormalized-energy configurations \cite{leble2017large,leble2017local,rougerie2016higher}. There are also interesting rigidity and tolerance properties that emerge in the $N \to \infty$ limit \cite{GP17, DV21, DHLM21, thoma2023nonrigidity}.
	
	This paper will be about an interacting particle system at high temperature, with a general pair-wise interaction. The study of interacting particle systems with a general interaction is a classical subject in statistical mechanics \cite{georgii2011gibbs,ruelle1967variational,ruelle1968statistical}, and also an active field of study \cite{garcia2019large, lambert2021poisson, chafai2014first,garcia2022generalized}. Particle systems at high temperature have also recently drawn attention, and this subject occasionally overlaps with the analysis of general interactions \cite{akemann2019high,chafai2014first,garcia2019large,hardy2021clt,padilla2020large,padilla2022large}.
	
	The analysis of interacting particle systems with a general pair-wise interaction is also linked to AI and machine learning, more specifically to neural networks. Neural networks can be used to accurately represent high-dimensional functions: given a high-dimensional function $f$, it may be represented as
	\begin{equation}
		f (x) := \lim_{n \to \infty} \frac{1}{n} \sum_{i=1}^{n} {\varphi}_{i} (x, \theta_{i}),
	\end{equation}
	where ${\varphi}_{i} (x, \theta_{i})$ are given functions, depending on the parameter $\theta_{i}$. One of the most frequent algorithms used to determine the parameters $\theta_{i}$ is stochastic gradient descent (SGD). Despite its ubiquity, very few rigorous results for convergence existed for SGD until recently. The approach in \cite{rotskoff2018neural,rotskoff2018parameters,rotskoff2022trainability,wen2024coupling} is to model the evolution of the parameters $\theta_{i}$ under the SGD as a particle system with an evolution given by an SPDE. By modelling parameters as a particle system, the energy landscape for the empirical measure becomes convex, and the authors are able to show a convergence rate of $O(n^{-1})$ to the mean-field limit. In this model, the interaction between the particles depends on the error between the measurements and the approximating function, and on the functions ${\varphi}_{i}$. In general, however, it is not a Coulomb or Riesz interaction. If the ${\varphi}_{i}$ are radial basis function networks, the interaction between the particles may be given by a Gaussian kernel, which is weakly interacting (see Definition \ref{def:weakint}). Given that the function $f$ is high dimensional, this approach is linked to particle systems for large $d$ as well. 
	
	The main result in this paper is that, for a specific class of interaction kernels, which we call weakly interacting, the local point process converges to a Poisson process for $ N^{-1} \ll \beta \ll N^{-\frac{1}{2}}$, a wider temperature than the one considered in \cite{lambert2021poisson} (which, to our knowledge is currently the most general result in that direction), namely $\beta \simeq N^{-1}$. Hence, we prove that the hypothesis on the temperature scaling may be weakened for this specific class of interaction kernels. Our ``weakly interacting" kernels are essentially kernels with an integrable positive Fourier transform, see Definition \ref{def:weakint}. 
	
	\section{Main results} \label{s:2}
	
	We will now state precisely the main results of the paper. We begin by defining exactly the class of confining potentials and interaction kernels that we deal with. 
	
	\begin{definition}
		\label{def:weakint}
		An interaction kernel $g$ is called weakly interacting if $g$ satisfies:
		\begin{itemize}
			\item[1.] $g(x)=g(-x)$.
			\item[2.] $\widehat{g} > 0$ a.e.
			\item[3.] $\widehat{g} \in L^{1}$.
		\end{itemize}
		Note that, as a consequence of item $3$, $g$ is bounded below (and above). Without loss of generality, we assume that this lower bound is $0$. Also as a consequence of item $3$, $g$ is continuous. Hence, we may assume without loss of generality the following:
		\begin{itemize}
			\item[4.] $g >0$ a.e. 
			\item[5.] $g$ is continuous.
		\end{itemize}
		
		A confining potential $V$ is called admissible if $V$ satisfies:
		\begin{itemize}
			\item[1.] $V$ is l.s.c.
			\item[2.] $\lim_{x \to \infty} V(x) = \infty$.
            \end{itemize}
            Note that as a consequence of items 1 and 2, $V$ satisfies the following additional properties:
            \begin{itemize}
			\item[3.] $\int_{|x|\geq 1}e^{-\alpha V(x)}\, \mathrm d x<+\infty$ for large enough $\alpha \geq \alpha_0$. 
            \item[4.] $V$ is bounded below. W.L.O.G. we assume that this bound is $0$, i.e. That $V$ is non-negative.
		\end{itemize}
            
	\end{definition}
	
	\begin{definition}
		A key quantity throughout the paper will be the potential field $h^\mu := g \ast \mu$ associated to a finite measure $\mu$. We often take $\mu = \emp_N$.
	\end{definition}
	We will now state our results for weakly interacting kernels and admissible confining potentials. As mentioned in the introduction, the main goal of the paper is to prove that the local point process converges to a Poisson point process. Along the way, we prove concentration bounds for the potential field of $\emp_N$ at given points, and asymptotes for the Laplace transform of fluctuations. These are interesting results in their own right, and so we state them in this section.
	
	We recall the definition of the equilibrium measure $\mu_V$ in (\ref{eq:defmuV}) and thermal equilibrium measure $\mu_\theta$ in (\ref{eq:defmutheta}). Throughout the paper, we assume $\theta \geq \alpha_0$ from Definition \ref{def:weakint}.
	
	\begin{proposition}[Concentration bounds for the potential field]
		\label{prop:concentration}
		Assume that $g$ is weakly interacting, and $V$ is admissible. For any $k \in \mathbf{N}$, $\epsilon > 0$, and points $y_{1}, y_{2}, ... y_{k} \in \mathbf{R}^{d}$, we have
		\begin{equation}
			\mathbb{P}_{N, \beta} \left( \left|  \sum_{i=1}^{k} h^{\rm emp_{N} - \mu_{\theta}} (y_{i}) \right| >  \epsilon\right) \leq \exp\(N \beta g(0) - \frac{N^2 \beta \epsilon^2}{g(0) k^2} \).
		\end{equation}
	\end{proposition}
	The proof is found in Section \ref{sect:proofconc}., and as a consequence we derive in Section \ref{sect:prooflaplace} the following bound.

\begin{proposition}[Asymptotics of the Laplace transform of field fluctuations] 
	\label{prop:laplace}
	Assume that $g$ is weakly interacting, and $V$ is admissible. Let $k \in \mathbf{N}$ and $y_{1}, y_{2}, ... y_{k} \in \mathbf{R}^{d}$. Then 
	\begin{equation}
		\mathbb{E}_{N, \beta} \left[ \exp\left( -N \beta \left( \sum_{i=1}^{k} h^{\rm emp_{N}} (y_{i}) \right) \right) \right] = M_{N} \exp\left( -N \beta \left( \sum_{i=1}^{k} h^{\mu_{\theta}} (y_{i}) \right) \right) + A_{N}, 
	\end{equation}
	where $M_N > 0$ and $A_N$ satisfy
	\begin{equation}
		\begin{split}
			|\log M_{N}| &\leq {\sqrt{2N} g(0) k}, \\
			|A_{N}| &\leq \exp(-N \beta g(0))\left(1 + e^{\sqrt{2N} g(0) k} \right).
		\end{split}
	\end{equation}
\end{proposition}
	
	We now state the main result of this paper: that the local point process converges to a Poisson point process. 
	\begin{theorem}[Convergence to Poisson point process]
		\label{teo:poisson}
		Assume that $g$ is weakly interacting, and $V$ is admissible. Assume also that $\mathrm{ent}[\mu_{V}] < \infty$. Let $x^{*} \in \mathbf{R}^{d}$ and define the local point process $\Xi$ by 
		\begin{equation}
			\Xi := \sum_{i=1}^{N} \delta_{N^{\frac{1}{d}} \left( x_{i} - x^{*} \right)}.
		\end{equation}
		
		Then, for a.e. $x^{*}$, if $N^{-1} \ll \beta \ll N^{-\frac{1}{2}}$, $\Xi$ converges to a Poisson point process of intensity $\mu_{V}(x^{*})$ as $N \to \infty$.  
	\end{theorem}
	The proof is found in Section \ref{sec:mainteo}. 
	
	Our proof also gives a precise asymptotic for the first marginal of the Gibbs measure, defined as
	\begin{equation}
		\rho(x) :=  \frac{1}{Z^{V}_{N, \beta}} \int_{\mathbf{R}^{d-1}} \exp \left( - \beta \mathcal{H}_{N}(x, X_{N-1}) \right) \, \mathrm d X_{N-1}. 
	\end{equation}
	The problem of finding an asymptotic for the first marginal was also addressed in \cite{rougerie2014quantum}, motivated by a link to the  fractional quantum Hall effect.  
	
	\begin{corollary}[First marginal and confinement]
		
		We have
		\begin{equation}
			\rho(x) = M''_{N}\mu_{\theta}(x) + A''_{N},
		\end{equation}
		for some $M'_{N}, A'_{N}$ satisfying 
		\begin{equation}
			\label{eq:marginal}
			\begin{split}
				|\log M'_{N}| &\leq C\beta \sqrt{N} \\
				|A'_{N}| &\leq C\exp \left( - C^{-1} \beta N \right )
			\end{split}
		\end{equation}
		for some constant $C > 0$ depending only on $V$ and $g$.
		
		If $V$ grows polynomially fast at $\infty$, and if $\beta = N^{-s }$ for some $s \in \left( \frac{1}{2}, 1 \right)$, then equation \eqref{eq:marginal} implies a confinement bound, or a bound on the probability that there is one particle outside of a compact set $K$. To be precise, let $K \subset \mathbf{R}^{d}$ be a compact set such that $V(x) + h^{\mu_{V}}(x)- c_{\infty}$ (see Lemma \ref{lem:equilibrium}) is bounded below by a positive constant outside of $K$. Then
		\begin{equation}
			\begin{split}
				\mathbb{P}_{N, \beta} \left( \emp_{N}\left( \mathbf{R}^{d} \setminus K \right) \geq \frac{1}{N} \right) &\leq N \int_{\mathbf{R}^{d} \setminus K} \rho(x) \, \mathrm d x\\
				&\leq \exp \left( - c N^{\gamma} \right),
			\end{split}    
		\end{equation}
		where $c,\gamma >0$ depend on $V,g,$ and $K$.
	\end{corollary}
	
	In the Coulomb case, results from \cite{armstrong2022thermal} imply that for any compact set $K$ that contains a neighborhood of the support of $\mu_{V}$, we have that $V(x) + h^{\mu_{V}}(x)- c_{\infty}$ is bounded below by a positive constant outside of $K$.
	
	\begin{remark}
		Physical considerations lead one to expect that the optimal regime for convergence of the local point process to a Poisson process is $\beta \to 0$ as $N \to \infty$ at any speed. Proving this in a mathematically rigorous way is an interesting direction for future work. 
	\end{remark}
	
	\begin{remark}
		Our analysis does not cover Coulomb, Riesz, or $1d-$log gases, since their Fourier transform is not integrable. Proving that the local point process of a Coulomb, Riesz, or $1d-$log gas converges to a Poisson process in a wider temperature regime than $\beta \simeq \frac{1}{N}$ is work in progress. 
	\end{remark}
	
	\section{Preliminaries}
	\label{sect:prelims}
	
	Before giving the proof of the main results, we establish some foundational results in interacting particle systems. Some of these results may exist in the literature, but since we were unable to find a source corresponding to our specific case we present arguments for them here. To our knowledge, Lemma \ref{lem:convergence} and Remark \ref{rem:qual} are new.
	
	For the rest of the paper, we will commit the abuse of notation of not distinguishing between a measure and its density. 
	
	\begin{lemma}[Equilibrium measure]
		\label{lem:equilibrium}
		Assume that $g$ is weakly interacting, and $V$ is admissible. Then the functional 
		\begin{equation}
			\mathcal{E}_{V} (\mu) := \mathcal{E}(\mu) + \int_{\mathbf{R}^{d}} V \, \mathrm d \mu 
		\end{equation}
		has a unique minimizer in the space of probability measures which we call the \emph{equilibrium measure} and denote $\mu_{V}$. The equilibrium measure has compact support, denoted $\Sigma$, and satisfies the Euler-Lagrange equation
		\begin{equation}
			\label{ELeq}
			\begin{split}
				h^{\mu_{V}} + \frac{V}{2} - c_{\infty} \geq &0 \text{ in } \mathbf{R}^{d}\\
				h^{\mu_{V}} + \frac{V}{2} - c_{\infty} = &0 \:\mu_V-\text{a.e. in } \{x \in \mathbf{R}^{d} | \mu_{V}(x)>0\}, 
			\end{split}
		\end{equation}
		where 
		\begin{equation}
			\label{eq:equilibriumconstant}
			c_{\infty} = \mathcal{E}(\mu_{V}) + \frac{1}{2} \int_{\mathbf{R}^{d}} V \, \mathrm d \mu.
		\end{equation}
		
	\end{lemma}
	The proof is a standard argument, and versions can be found for the two dimensional log-gas in \cite[Chapter 1]{saff97pot}, for general interactions on compact sets in \cite[Chapter 4]{borodachov19pot} and for a general Riesz gas in \cite[Chapter 2]{serfaty15gl}.

			\begin{lemma}[Thermal equilibrium measure]
				
				Assume that $g$ is weakly interacting, $V$ is admissible, and let $\theta>0$. Then the functional 
				\begin{equation}
					\mathcal{E}_{\theta} (\mu) := \mathcal{E}(\mu) + \int_{\mathbf{R}^{d}} V \, \mathrm d \mu + \frac{1}{\theta} \mathrm{ent}[\mu]
				\end{equation}
				has a unique minimizer in the space of probability measures which we call \emph{thermal equilibrium measure} and denote $\mu_{\theta}$. The thermal equilibrium measure is everywhere positive, and satisfies the Euler-Lagrange equation
				\begin{equation}\label{foc}
					\mu_{\theta} = L_{\theta}^{-1} \exp \left( - \theta \left( 2h^{\mu_{\theta}} +V \right) \right),
				\end{equation}
				where
				\begin{equation}
					\label{eq:thermalconstant}
					L_{\theta} = \int_{\mathbf{R}^{d}} \exp \left( - \theta \left( 2h^{\mu_{\theta}}(x) +V(x) \right) \right) \, \mathrm d x. 
				\end{equation}    
			\end{lemma}
			\begin{proof}
				The existence of the thermal equilibrium measure follows for instance from \cite[Lemma 2.1]{armstrong2022thermal}, which only uses that $g$ is bounded from below and that $V \rightarrow +\infty$ and satisfies the exponential integrability assumption
				\begin{equation}
					\int_{|x|\geq 1}e^{-\frac{\theta}{2}V(x)}\, \mathrm d x<+\infty.
				\end{equation}
				(\ref{foc}) then follows immediately from the form of the thermal equilibrium measure given in \cite[Proposition A.2]{lambert2021poisson}, since our kernel $g$ satisfies all of the hypotheses of that proposition.
			\end{proof}
			Once we have existence of the thermal equilibrium measure, $N^2\mathcal{E}_{\theta}$ gives the leading order asymptotics of the Hamiltonian. To analyze next-order behavior, it is then useful to split off this deterministic term and consider the remainder. This can be accomplished by the following \textit{splitting formula} and its corollary; the version stated below is exactly \cite[Lemma 2.1]{armstrong2021local}, which only relies on the form of the Hamiltonian $\mathcal{H}_N$.
			\begin{lemma}[Thermal splitting formula]
				\label{lem:splitting}
				Let $\theta >0$, and define
				\begin{equation}
					\zeta_{\theta} := - \frac{1}{\theta} \log(\mu_{\theta}).
				\end{equation}
				
				Then for any point configuration $X_{N} \in \mathbf{R}^{d \times N}$ the Hamiltonian $\mathcal{H}_{N}$ can be rewritten (split) as  
				\begin{equation}\label{eq:thermspltfrm}
					\mathcal{H}_{N} (X_{N}) = N^{2} \left(  \mathcal{E}_{\theta} (\mu_{\theta}) +  \F_{N}(X_{N}, \mu_{\theta})+ N\sum_{i=1}^N \zeta_{\theta}(x_i) \right),
				\end{equation}
				where, given a measure $\mu$,
				\begin{equation}
					\F_{N}(X_{N}, \mu) := \frac{1}{N^{2}} \sum_{i \neq j} g(x_{i} - x_{j}) + \mathcal{E}(\mu) - \frac{2}{N} \sum_{i=1}^{N} h^{\mu}(x_{i}). 
				\end{equation}
			\end{lemma}
			A computation then immediately yields the following equivalent definition of the Gibbs measure.
			
			\begin{corollary} 
				The Gibbs measure may be rewritten as 
				\begin{equation}
					\mathrm d \mathbb{P}_{N, \beta} (X_{N}) = \frac{1}{K_{N, \beta}^{\theta}} \exp \left( - N^{2} \beta \left[ \F_{N}(X_{N}, \mu_{\theta}) \right]  \right) \prod_{i=1}^{N} \mu_{\theta}(x_{i}) \, \mathrm d X_{N},
				\end{equation}
				where
				\begin{equation}
					K_{N, \beta}^{\theta} = \frac{ Z^{V}_{N, \beta} }{\exp \left( N^{2} \beta \mathcal{E}_{\theta}(\mu_{\theta}) \right)}.
				\end{equation}    
				
				We call $K_{N, \beta}^{\theta}$ the \emph{next-order partition function}. A useful fact is that its logarithm is always non-negative.
			\end{corollary}
			
			\begin{lemma} \label{lem:logK}
				Assume that $g$ is weakly interacting, and $V$ is admissible. Then the next-order partition function is greater than $1$, i.e.
				\begin{equation}
					\log K_{N, \beta}^{\theta} >0.
				\end{equation}
			\end{lemma}
			\begin{proof}
				The proof is as in \cite[Proposition 5.10]{padilla2023concentration}, which considered the analogous question for the Coulomb gas. The result follows immediately from the definition of the Gibbs measure as the probability measure that minimizes the free energy functional 
				\begin{equation}
					\mathbb{Q} \mapsto \int_{\mathbf{R}^{d \times N}}\mathcal{H}_N(X_{N})\, \mathrm d \mathbb{Q}\left(X_N\right)+\frac{1}{\beta}\mathrm{ent}[\mathbb{Q}]
				\end{equation}
				over all $\mathbb{Q} \in \mathcal{P}\left(\mathbf{R}^{d \times N}\right)$. Inserting $\mathbb{Q}=\mathbb{P}_{N,\beta}$ yields
				\begin{equation}
					-\frac{\log Z_{N,\beta}}{\beta}=\inf_{\mu \in\mathcal{P}\left(\mathbf{R}^{d \times N}\right) }\int_{\mathbf{R}^{d \times N}}\mathcal{H}_N(X_{N})\, \mathrm d \mu\left(X_N\right)+\frac{1}{\beta}\mathrm{ent}[\mu].
				\end{equation}
				Inserting instead $\mu_{\theta}^{\otimes N}$ and computing yields
				\begin{equation}
					-\frac{\log Z_{N,\beta}}{\beta}\leq  N^2 \mathcal{E}_{\theta}(\mu_{\theta})-\frac{N}{2}\mathcal{E}(\mu_{\theta}) <N^2 \mathcal{E}_\theta(\mu_{\theta})
				\end{equation}
				by positivity of $g$. The splitting formula (\ref{eq:thermspltfrm}) then yields
				\begin{equation}
					-\frac{\log Z_{N,\beta}}{\beta}=\frac{-\log K_{N,\beta}^\theta+\beta N^2 \mathcal{E}(\mu_{\theta})}{\beta} \leq N^2 \mathcal{E}_\theta(\mu_{\theta}) \implies \log K_{N,\beta}^\theta> 0.
				\end{equation}
			\end{proof}
			
			From this point forward, we present results on the thermal equilibrium measure that are not essentially taken from existing literature. 
			
			\begin{lemma}[Convergence]
				\label{lem:convergence}
				Assume that $g$ is weakly interacting, and $V$ is admissible. Suppose also that $ \mathrm{ent}[\mu_{V}] < \infty$. Then, as $\theta \to \infty$, $\mu_{\theta}$ converges to $\mu_{V}$ in the $L^{1}$ topology, and $h^{\mu_{\theta}}$ converges to $h^{\mu_{V}}$ in the $L^{\infty}$ topology. Furthermore, for a.e. $x \in \mathbf{R}^{d}$ and any sequence $\delta(\theta)$ such that $\lim_{\theta \to \infty } \delta = 0$,
				\begin{equation}
                \label{eq:microconv}
					\lim_{\theta \to \infty} \fint_{B(x, \delta)} \left| \mu_{\theta}(s) - \mu_{V}(x) \right| \, \mathrm d s  = 0.
				\end{equation}
			\end{lemma}
			\begin{proof}
				In the case where $g$ is the Coulomb interaction, very precise convergence is established by connections with the classical obstacle problem (cf. \cite[Theorem 1]{armstrong2022thermal}). However, that machinery is not available to us in the general interaction case.
				
				\textbf{Step 1:} We will first prove that $\mu_{\theta} \to \mu_{V}$ in the topology of weak convergence of probability measures. We will use the language of $\Gamma$-convergence (see \cite{braides2002gamma}). In particular, a sequence of functionals $F_n:X \rightarrow \mathbb{R}$ on a metric space is said to $\Gamma$-\textit{converge} to a functional $F_\infty$ if 
				\begin{enumerate}
					\item ($\Gamma-\liminf$ inequality) For all sequences $x_n \rightarrow x_\infty$, 
					\begin{equation*}
						F_\infty(x_\infty)\leq \liminf_{n \rightarrow \infty}F_n(x_n)
					\end{equation*}
					\item (recovery sequence) There exists a sequence $x_n \rightarrow x_\infty$ such that 
					\begin{equation*}
						F_\infty(x_\infty)\geq \limsup_{n \rightarrow \infty}F_n(x_n)    \end{equation*}
				\end{enumerate}
				The reason for using this framework is that $\Gamma$-convergence implies weak convergence of the minimizers, as in \cite[Theorem 1.21]{braides2002gamma} and \cite[Theorem 2.2]{serfaty15gl}.
				
				\textbf{Step 1.1:} $\Gamma-\liminf$ inequality. 
				
				The lack of positivity of the entropy presents a difficulty, which we circumvent by a trick as in \cite{armstrong2022thermal}; namely, we can rewrite the function $\mathcal{E}_\theta$ as
				\begin{equation}
					\mathcal{E}_\theta(\mu)=\mathcal{E}(\mu)+\left(1-\frac{\alpha_0}{\theta}\right)\int_{\mathbf{R}^{d}} V\, \mathrm d \mu+\left(\int \frac{\alpha_0}{\theta}V\, \mathrm d \mu +\mathrm{ent}[\mu]\right),
				\end{equation}
				where $\alpha_0$ is as in Definition \ref{def:weakint}. 
				Consider a sequence $\theta_n \rightarrow +\infty$ and suppose that $\mu_n \rightharpoonup \mu$. Then, since $g$ and $V$ are both lower semicontinuous and bounded from below, we have for any $\epsilon>0$ that 
				\begin{align}
					\liminf_{n \rightarrow \infty}\left( \mathcal{E}(\mu_{n})+\int_{\mathbf{R}^{d}} \left(1-\frac{\alpha_0}{\theta_n}\right)V\, \mathrm d \mu_n \right) &\geq\liminf_{n \rightarrow \infty}\left( \mathcal{E}(\mu_{n})+\int_{\mathbf{R}^{d}} \left(1-\epsilon\right)V\, \mathrm d \mu_n \right) \\ \notag
					&\geq  \mathcal{E}(\mu)+\int_{\mathbf{R}^{d}} (1-\epsilon)V~ d\mu.
				\end{align}   
				Now, if $\mu_n$ is not absolutely continuous with respect to Lebesgue measure for all $n \geq N$ for some $N$, then $\mathrm{ent}[\mu_n]$ is infinite for all $n \geq N$ and so 
				\begin{equation}
					\liminf_{n \rightarrow \infty}\left(\int_{\mathbf{R}^{d}} \frac{\alpha_0}{\theta_n}V\, \mathrm d \mu_n +\frac{1}{\theta_n}\mathrm{ent}[\mu_n]\right)=+\infty \geq 0.
				\end{equation}
				Otherwise, extracting the subsequence $\mu_{n_k}$ that are absolutely continuous to Lebesgue measure yields
				\begin{align*}
					\liminf_{n \rightarrow \infty}\left(\int_{\mathbf{R}^{d}} \frac{\alpha_0}{\theta_n}V\, \mathrm d \mu_n +\frac{1}{\theta_n}\mathrm{ent}[\mu_n]\right) &\geq \liminf_{k \rightarrow \infty} \int_{\R^d}\left(\frac{\alpha_0 V}{\theta_{n_k}}+\frac{1}{\theta_{n_k}}\log \mu_{n_k}\right)~d\mu_{n_k} \\
					&\geq \liminf_{k \rightarrow \infty}\frac{-1}{\theta_{n_k}}\int_{\R^d}e^{-\alpha_0V(x)-1}~dx=0
				\end{align*}
				where we have used the inequality
				\begin{equation*}
					\frac{c\gamma}{\theta}+\frac{\gamma \log \gamma}{\theta} \geq -\frac{1}{\theta}e^{-c-1},
				\end{equation*} 
				with $c=\alpha_0 V$ and $\gamma=\mu_{n_k}$,
				which can be seen by minimizing the above as a function of $\gamma$. 
				
				In either case,
				\begin{equation*}
					\liminf_{n \rightarrow \infty}\mathcal{E}_{\theta_n}(\mu_n)\geq  \mathcal{E}(\mu)+\int_{\mathbf{R}^{d}} (1-\epsilon)V~ d\mu.
				\end{equation*}
				
				Since $\epsilon>0$ was arbitrary, we conclude by monotone convergence that $\liminf_{n \rightarrow +\infty}\mathcal{E}_{\theta_n}(\mu_n)\geq \mathcal{E}_V(\mu)$. 
				
				\textbf{Step 1.2:} $\Gamma-\limsup$ inequality. 
				
				Under the hypothesis that $ \mathrm{ent}[\mu_{V}] < \infty$, it is trivially true that
				\begin{equation}
					\lim_{\theta \to \infty} \mathcal{E}_{\theta}(\mu_{V}) = \mathcal{E}_{V}(\mu_{V}). 
				\end{equation}
				
				\textbf{Step 1.3:} Uniform coercivity and conclusion of step 1.  
				
				Now we show that $\{\mu_{\theta_n}\}$ is a tight sequence of probability measures. This will follow from the boundedness of the sequence $\{\mathcal{E}_{\theta_n}(\mu_{\theta_n})\}$ in the same way as Lemma 4.1 (following \cite[Lemma 2.10]{serfaty15gl}). Without loss of generality, assume $\theta_1=\inf \theta_n$. First, the boundedness follows from $\frac{1}{\theta_n}\rightarrow 0$ and the minimality of $\mu_{\theta_n}$:
				\begin{equation}
					\mathcal{E}_{\theta_n}(\mu_{\theta_n}) \leq \mathcal{E}_{\theta_n}(\mu_{\theta_1}) \leq \mathcal{E}_{\theta_1}(\mu_{\theta_1}):=C_1<+\infty
				\end{equation}
				for all $n$.  Now, for any $C_2>0$, there is a compact set $K \times K$ outside of which 
				\begin{equation}
					g(x-y)+\gamma_n\frac{V(x)}{2}+\gamma_n\frac{V(y)}{2}>C_2
				\end{equation}
				with $\gamma_n=\left(1-\frac{\alpha_0}{\theta_n}\right)$
				independently of $n$ by our growth assumption on $V$ (we can choose this independently of $n$ since $\theta_n \rightarrow +\infty$). Then, moving some of the potential onto the entropy as above, we have
				\begin{align*}
					C_1 &\geq \mathcal{E}(\mu_{\theta_n})+\left(1-\frac{\alpha_0}{\theta_n}\right)\int_{\mathbf{R}^{d}} V\, \mathrm d \mu_{\theta_n}+\int_{\mathbf{R}^{d}} \left(\frac{\alpha_0}{\theta_n}V+\frac{1}{\theta_n}\log \mu_{\theta_n}\right)\, \mathrm d \mu_{\theta_n} \\
					&\geq \iint_{\R^d \times \R^d} \left(g(x-y)+\gamma_n\frac{V(x)}{2}+\gamma_n\frac{V(y)}{2}\right)\, \mathrm d \mu_{\theta_n}(x)\, \mathrm d \mu_{\theta_n}(y)-\frac{1}{\theta_n}\int_{\mathbf{R}^{d}} e^{-\alpha_0V(x)-1}\, \mathrm d x \\
					& \geq -C_3+C_2\mu_{\theta_n}\otimes \mu_{\theta_n}((K \times K)^c) \\
					& \geq -C_3+C_2\mu_{\theta_n}(K^c),
				\end{align*}
				where we have used that $g(x-y)+\gamma_n\frac{V(x)}{2}+\gamma_n\frac{V(y)}{2}$ is everywhere bounded below  independently of $n$ since $|\gamma_n-1|$ is small, and $-\frac{1}{\theta_n}\int_{\mathbf{R}^{d}} e^{-\alpha_0V(x)-1}\, \mathrm d x$ is bounded below as well since $\theta_n \rightarrow \infty$. Since $C_2$ can be made arbitrarily large, this gives us the tightness of the sequence $\{\mu_{\theta_n}\}$.
				
				We can now conclude as in \cite[Theorem 1.21]{braides2002gamma}, or as in \cite[Theorem 2.2]{serfaty15gl}. Since $\{\mu_{\theta_n}\}$ is a tight sequence, it has a weak limit $\mu$. Furthermore, since $\mu_V$ has a density, we have
				\begin{equation}
					\mathcal{E}_V(\mu_V) = \limsup_{n \rightarrow +\infty}\mathcal{E}_{\theta_n}(\mu_V) \geq \limsup_{n \rightarrow +\infty}\mathcal{E}_{\theta_n}(\mu_{\theta_n}) \geq \mathcal{E}_V(\mu). 
				\end{equation}
				The second inequality follows from the fact that $\mu_{\theta_n}$ minimizes $\mathcal{E}_{\theta_n}$, and the third inequality is a result of lower semicontinuity. Since the equilibrium measure is unique, we must have $\mu=\mu_V$ and we see that $\mu_{\theta_n}$ converges weakly to $\mu_V$ under the topology of weak convergence.
				
				\textbf{Step 2:} Now, we prove that $\mu_{\theta} \to \mu_{V}$ in the $L^{1}$ topology.
				
				Note that, by definition of $\mu_{\theta}$, we have that, for any $\theta >0$,
				\begin{equation}
					\mathrm{ent}[\mu_{\theta}] \leq \mathrm{ent}[\mu_{V}] < \infty.
				\end{equation}
				In particular, 
				\begin{equation}
					\overline{\lim}_{\theta \to \infty}\mathrm{ent}[\mu_{\theta}] \leq \mathrm{ent}[\mu_{V}].
				\end{equation}
				On the other hand, by l.s.c.,
				\begin{equation}
					\mathrm{ent}[\mu_{V}] \leq  \underline{\lim}_{\theta \to \infty}\mathrm{ent}[\mu_{\theta}],
				\end{equation}
				which implies that
				\begin{equation}
					\label{eq:convent}
					{\lim}_{\theta \to \infty}\mathrm{ent}[\mu_{\theta}] = \mathrm{ent}[\mu_{V}].
				\end{equation}
				
				Given $\epsilon>0$, we define the probability measure $\mu_{V}^{\epsilon}$ as 
				\begin{equation}
					\mu_{V}^{\epsilon} := \frac{\mu_{V} + \epsilon \exp\left( - \alpha_{0} V \right) }{\int_{\mathbf{R}^{d}} \mu_{V} + \epsilon  \exp\left( - \alpha_{0} V \right) },
				\end{equation}
				where $\alpha_0$ is as in Definition \ref{def:weakint}. Then, by Pinsker's inequality,
				\begin{equation}
					\begin{split}
						\limsup_{\theta \to \infty} \left\| \mu_{\theta} -  \mu_{V}^{\epsilon} \right\|_{L^{1}} & \leq \limsup_{\theta \to \infty}   \sqrt{\int_{\mathbf{R}^{d}} \log \left( \frac{\mu_{\theta}}{\mu_{V}^{\epsilon}} \right) \mu_{\theta} }\\
						& = \limsup_{\theta \to \infty}   \sqrt{\mathrm{ent}[\mu_{\theta}] - \int_{\mathbf{R}^{d}} \log \left( {\mu_{V}^{\epsilon}} \right) \mu_{\theta} }\\
						& =  \sqrt{\mathrm{ent}[\mu_{V}] - \int_{\mathbf{R}^{d}} \log \left( {\mu_{V}^{\epsilon}} \right) \mu_{V} },
					\end{split}
				\end{equation}
				where we've used the weak convergence of $\mu_\theta$ to $\mu_V$ established above to conclude that $\int \log(\mu_V^\epsilon)\mu_\theta \rightarrow \int \log(\mu_V^\epsilon)\mu_V$.
				
				Furthermore, by triangle inequality, there holds 
				\begin{equation}
					\limsup_{\theta \to \infty} \left\| \mu_{\theta} -  \mu_{V} \right\|_{L^{1}} \leq \sqrt{\mathrm{ent}[\mu_{V}] - \int_{\mathbf{R}^{d}} \log \left( {\mu_{V}^{\epsilon}} \right) \mu_{V} } + \left\| \mu_{V} -  \mu_{V}^{\epsilon} \right\|_{L^{1}}.
				\end{equation}
				Letting $\epsilon$ tend to $0$, we may conclude that $\mu_{\theta} \to \mu_{V}$ in the $L^{1}$ topology.  
				
				Note that $\left\| g \right\|_{L^{\infty}} < \infty$ by property 3 of Definition \ref{def:weakint}. Hence, Young's convolution inequality, implies that $h^{\mu_{\theta}} \to h^{\mu_{V}}$ in $L^{\infty}$. 
				
				\textbf{Step 3:} Finally, we prove that for a.e. $x \in \mathbf{R}^{d}$ and any sequence $\delta(\theta)$ such that $\lim_{\theta \to \infty } \delta = 0$,
				\begin{equation}
					\lim_{\theta \to \infty} \fint_{B(x, \delta)} \left| \mu_{\theta}(s) - \mu_{V}(x) \right| \, \mathrm d s  = 0.
				\end{equation}
				
				To prove this claim, note that, for any $\theta>0$ and a.e. $x$, we have
				\begin{equation}
					\begin{aligned}
						\lefteqn{ \fint_{B(x, \delta)} \left| \mu_{\theta}(s) - \mu_{V}(x) \right| \, \mathrm d s } \quad & \\ & 
						\leq  \fint_{B(x, \delta)}  \left|\mu_{\theta}(s)  - \mu_{V}(s) \right| \, \mathrm d s  + \fint_{B(x, \delta)}  \left| \mu_{V}(s)  - \mu_{V}(x) \right| \, \mathrm d s \\
						& \leq  M(\mu_{\theta} - \mu_{V})(x) +\fint_{B(x, \delta)}  \left| \mu_{V}(s)  - \mu_{V}(x) \right| \, \mathrm d s.
					\end{aligned}
				\end{equation}
				
				where $M(\cdot)$ denotes the Hardy-Littlewood maximal function. 
				
				Let $\epsilon > 0$. Then 
				\begin{equation}
					\begin{aligned}
						\lefteqn{ \overline{\lim}_{\theta \to \infty } \left| \left\{ x \in \mathbf{R}^{d} : \fint_{B(x, \delta)} \left| \mu_{\theta}(s) - \mu_{V}(x) \right| \, \mathrm d s  > \epsilon \right\} \right|} \quad & \\ &
						\leq  \overline{\lim}_{\theta \to \infty } \left| \left\{ x \in \mathbf{R}^{d} : \left| M(\mu_{\theta} - \mu_{V})(x) \right| > \epsilon \right\} \right| \\
						& \quad + \overline{\lim}_{\theta \to \infty } \left| \left\{ x \in \mathbf{R}^{d} : \fint_{B(x, \delta)}  \left| \mu_{V}(s)  - \mu_{V}(x) \right| \, \mathrm d s > \epsilon \right\} \right|.
					\end{aligned}    
				\end{equation}
				
				By the Hardy-Littlewood maximal inequality, and step 2, there holds
				\begin{equation}
					\begin{split}
						\overline{\lim}_{\theta \to \infty} \left| \left\{ x \in \mathbf{R}^{d} : \left| M(\mu_{\theta} - \mu_{V})(x) \right| > \epsilon \right\} \right| &\leq  \overline{\lim}_{\theta \to \infty} \frac{C_{d}}{\epsilon} \left\| \mu_{\theta} - \mu_{V} \right\|_{L^{1}} = 0.
					\end{split}
				\end{equation}
				
				On the other hand, Lebesgue's Theorem implies that 
				\begin{equation}
					\overline{\lim}_{\theta \to \infty } \left| \left\{ x \in \mathbf{R}^{d} : \fint_{B(x, \delta)}  \left| \mu_{V}(s)  - \mu_{V}(x) \right| \, \mathrm d s > \epsilon \right\} \right| =0.
				\end{equation}
				Since $\epsilon$ is arbitrary, we may conclude. 
				
			\end{proof}
			
			\begin{remark}[Qualitative behaviour]
				\label{rem:qual}
				The constant $L_{\theta}$ defined in \eqref{eq:thermalconstant} satisfies that
				\begin{equation}
					\label{eq:limLtheta}
					\lim_{\theta \to \infty} -\frac{1}{\theta} \log L_{\theta} = 2 c_{\infty},
				\end{equation}
				where $c_{\infty}$ is given by \eqref{eq:equilibriumconstant}. Furthermore, for all $\theta>0$ there exists a constant $c_{\theta}$ such that, for any $\theta'\geq \theta$ there holds
				\begin{equation}
					\mu_{\theta'} (x) \leq \exp \left( -\theta' \left( 2h^{\mu_{V}}(x) + V(x) - 2 c_{\theta}\right) \right).
				\end{equation}
				Also, $\lim_{\theta \to \infty} c_{\theta} = c_{\infty}$. 
			\end{remark}
			
			\begin{proof}
				Using Lemma \ref{lem:convergence}, and equations \eqref{ELeq} and \eqref{foc} we have that
				\begin{equation}
					\begin{split}
						\lim_{\theta \to \infty} -\frac{1}{\theta} \log L_{\theta} &= \inf_{x \in \mathbf{R}^{d}} 2 h^{\mu_{V}}(x) + V(x)\\
						&= 2c_{\infty}.
					\end{split}
				\end{equation}
				
				The rest of the remark follows from equations \eqref{foc}, \eqref{eq:limLtheta} by taking 
				\begin{equation}
					c_{\theta} := \sup_{\theta' \geq \theta}-\frac{1}{2 \theta'} \log L_{\theta'} + \left\|h^{\mu_{V}} - h^{\mu_{\theta'}}\right\|_{L^{\infty}}.
				\end{equation}
			\end{proof}
			
			\section{Proof of Proposition \ref{prop:concentration}}
			\label{sect:proofconc}
			
			We now prove Proposition \ref{prop:concentration}, restated here for convenience. 
			
			\begin{proposition}[Concentration bounds for the potential field]
				Assume that $g$ is weakly interacting, and $V$ is admissible. For any $k \in \mathbf{N}$, $\epsilon > 0$, and points $y_{1}, y_{2}, ... y_{k} \in \mathbf{R}^{d}$, we have
				\begin{equation}
					\mathbb{P}_{N, \beta} \left( \left|  \sum_{i=1}^{k} h^{\rm emp_{N} - \mu_{\theta}} (y_{i}) \right| >  \epsilon\right) \leq \exp\(N \beta g(0) - \frac{N^2 \beta \epsilon^2}{g(0) k^2} \).
				\end{equation}
			\end{proposition}
			
			\begin{proof}
				\textbf{Step 1:} Starting point.
				
				Using the splitting formula (Lemma \ref{lem:splitting}), we have that 
				\begin{equation}
					\label{eq:splitting}
					\begin{split}
						\lefteqn{ \mathbb{P}_{N, \beta} \left( \left|  \sum_{i=1}^{k} h^{\rm emp_{N} - \mu_{\theta}} (y_{i}) \right| >  \epsilon \right) } \quad & \\
						&= \frac{1}{K_{N, \beta}} \int_{\left\{ X_{N} : \left|  \sum_{i=1}^{k} h^{\rm emp_{N} - \mu_{\theta}} (y_{i}) \right| >  \epsilon \right\} } \exp \left( - N^{2} \beta \F_{N}(\emp_{N}, \mu_{\theta}) \right) \d \mu_{\theta}^{\otimes N}\\
						&\leq  \frac{1}{K_{N, \beta}} \sup_{\left\{ X_{N} : \left|  \sum_{i=1}^{k} h^{\rm emp_{N} - \mu_{\theta}} (y_{i}) \right| >  \epsilon \right\} } \exp \left( - N^{2} \beta \F_{N}(\emp_{N}, \mu_{\theta}) \right) \\
						&=  \frac{1}{K_{N, \beta}} \sup_{\left\{ X_{N} : \left|  \sum_{i=1}^{k} h^{\rm emp_{N} - \mu_{\theta}} (y_{i}) \right| >  \epsilon \right\} } \exp \left( - N^{2} \beta\left[ \mathcal{E}(\emp_{N}- \mu_{\theta}) -\frac{1}{N}g(0) \right]\right) \\
						&=  \frac{1}{K_{N, \beta}} \exp \left( - N^{2} \beta\left[ \Phi_{y_{1},...y_{k}}(\epsilon) -\frac{1}{N}g(0) \right]\right), 
					\end{split}
				\end{equation}
				where
				\begin{equation}
					\Phi_{y_{1},...y_{k}}(\epsilon) := \inf_{\left\{ \mu \in \mathcal{M}(\mathbf{R}^{d}) : \left|  \sum_{i=1}^{k} h^{\mu } (y_{i}) \right| >  \epsilon \right\} }  \mathcal{E}(\mu),
				\end{equation}
				and $\mathcal{M}(\mathbf{R}^{d})$ denotes the space of Radon measures. We will also use the obvious variants of $\Phi$ with different values of $k$.
				
				\textbf{Step 2:} Solving the variational problem.
				
				We claim that for any $y_{1}, ...y_{k} \in \mathbf{R}^{d}$,
				\begin{equation}
					\label{eq:coercive}
					\Phi_{y_{1},...y_{k}}(\epsilon) \geq \frac{\Lambda \epsilon^{2}}{k^{2}},
				\end{equation}
				for some $\Lambda>0$. 
				
				\textbf{Substep 2.1:} We start by proving that
				\begin{equation}
					\Lambda := \Phi_{0}(1) >0.
				\end{equation}
				
				To prove this claim, note that for any $\mu \in \mathcal{M}(\mathbf{R}^{d})$,
				\begin{equation}
					\begin{split}
						h^{\mu}(0) = \int_{\mathbf{R}^{d}} \widehat{h^{\mu}}(\xi) \, \mathrm d \xi
						= \int_{\mathbf{R}^{d}} \widehat{\mu}(\xi) \widehat{g}(\xi) \, \mathrm d \xi. 
					\end{split}    
				\end{equation}
				On the other hand, 
				\begin{equation}
					\mathcal{E}(\mu) = \int_{\mathbf{R}^{d}} \left|\widehat{\mu}(\xi) \right|^{2}\widehat{g}(\xi) \, \mathrm d \xi.
				\end{equation}
				By Cauchy-Schwartz, and since $\widehat{g}$ is positive a.e., we find
				\begin{equation}
					\left( \int_{\mathbf{R}^{d}} \widehat{\mu}(\xi) \widehat{g}(\xi) \, \mathrm d \xi\right)^{2} \leq \left( \int_{\mathbf{R}^{d}} \left|\widehat{\mu}(\xi) \right|^{2}\widehat{g}(\xi) \, \mathrm d \xi \right) \left( \int_{\mathbf{R}^{d}} \widehat{g}(\xi) \, \mathrm d \xi\right),
				\end{equation}
				and therefore for any $\mu \in \mathcal{M}(\mathbf{R}^{d})$ such that $h^{\mu}(0)=1$, there holds
				\begin{equation}
					\mathcal{E}(\mu) \geq \left( \int_{\mathbf{R}^{d}} \widehat{g}(\xi) \, \mathrm d \xi\right)^{-1} = \frac{1}{g(0)}, 
				\end{equation}
				which finishes the substep, and furthermore proves $\Lambda \geq (g(0))^{-1}$.
				
				\textbf{Substep 2.2:} We now prove the claim stated at the beginning of step 2. 
				
				Note that if $\left|  \sum_{i=1}^{k} h^{\mu } (y_{i}) \right| >  \epsilon$ then $\left| h^{\mu } (y_{i}) \right| >  \frac{\epsilon}{k}$ for some $i \in \{1,...k\}$ by a pigeonhole argument. Clearly, $\Phi_{y_i} = \Phi_0$ by translation invariance, and $\epsilon \mapsto \Phi_0(\epsilon)$ is quadratic since $|h^\mu(0)| \geq \epsilon$ iff $|h^{\mu/\epsilon}(0)| \geq 1$. It follows
				\begin{equation}
					\begin{split}
						\Phi_{y_{1},...y_{k}}(\epsilon) \geq \Phi_{0}\left(\frac{\epsilon}{k}\right) =  \frac{\Phi_0(1)\epsilon^2}{k^2}
						= \frac{\Lambda \epsilon^{2}}{k^{2}}.
					\end{split}
				\end{equation}
				By plugging equation \eqref{eq:coercive} into the last line of equation \eqref{eq:splitting}, and using Lemma \ref{lem:logK} to bound $K_{N, \beta} \geq 1$, we may conclude Proposition \ref{prop:concentration}. 
			\end{proof}
			
			\section{Proof of Proposition \ref{prop:laplace}}
			\label{sect:prooflaplace}
			
			In this section, we prove Proposition \ref{prop:laplace}, restated here for convenience. 
			
			\begin{proposition}[Asymptotics of the Laplace transform of field fluctuations] 
				Assume that $g$ is weakly interacting, and $V$ is admissible. Let $k \in \mathbf{N}$ and $y_{1}, y_{2}, ... y_{k} \in \mathbf{R}^{d}$. Then 
				\begin{equation}
					\mathbb{E}_{N, \beta} \left[ \exp\left( -N \beta \left( \sum_{i=1}^{k} h^{\rm emp_{N}} (y_{i}) \right) \right) \right] = M_{N} \exp\left( -N \beta \left( \sum_{i=1}^{k} h^{\mu_{\theta}} (y_{i}) \right) \right) + A_{N}, 
				\end{equation}
				where $M_N > 0$ and $A_N$ satisfy
				\begin{equation}
					\begin{split}
						|\log M_{N}| &\leq {\sqrt{2 N} g(0) k}\\
						|A_{N}| &\leq \exp(-N \beta g(0))\left(1 + e^{\sqrt{2N} g(0) k} \right).
					\end{split}
				\end{equation}
			\end{proposition}
			
			\begin{proof}
				Let $\epsilon = \frac{\sqrt{2} g(0) k}{\sqrt{N}}$, and define the event
				$$
				G = \left \{ \left \lvert \sum_{i=1}^k h^{\emp_N}(y_i) - h^{\mu_\theta}(y_i) \right \rvert > \epsilon \right \}.
				$$
				Since $g$ is nonnegative, and therefore $h^{\emp_N}$ is as well, we may use Proposition \ref{prop:concentration} to see
				\begin{align*}
					0 \leq \mathbb{E}_{N, \beta} \left[ \exp\left( -N \beta \left( \sum_{i=1}^{k} h^{\rm emp_{N}} (y_{i}) \right) \right) \1_{G} \right] &\leq \P_{N,\beta}(G) \\ &\leq \exp \left( -\frac{N^2 \beta \epsilon^2}{g(0) k^2} + N \beta g(0)\right) \\
					&= \exp \left( -N \beta g(0) \right).
				\end{align*}
			Defining $ M_N$ appropriately, we have
			\begin{align}
				\mathbb{E}_{N, \beta} \left[ \exp\left( -N \beta \left( \sum_{i=1}^{k} h^{\rm emp_{N}} (y_{i}) \right) \right) \1_{G^c} \right] &= M_{N} \exp\left( -N \beta \left( \sum_{i=1}^{k} h^{\mu_{\theta}} (y_{i}) \right) \right) \P_{N,\beta}(G^c)
			\end{align}
			with
			$
			| \log  M_N| \leq N \beta \epsilon
			$.
				Similarly, using $h^\mu_{\theta} \geq 0$, we may bound
				$$
				0 \leq M_{N} \exp\left( -N \beta \left( \sum_{i=1}^{k} h^{\mu_{\theta}} (y_{i}) \right) \right) \P_{N,\beta}(G) \leq e^{-N \beta (g(0) - \epsilon)}.
				$$
			Assembling the above yields the desired result.
			\end{proof}
			
			\section{Proof of Theorem \ref{teo:poisson}}
			\label{sec:mainteo}
			
			In this section, we prove the main result of this paper, namely Theorem \ref{teo:poisson}. We will use some background on point processes, which can be reviewed in the Appendix. The main result that we will use is the following lemma:
			
			\begin{lemma}[Weak convergence and correlation functions]
				\label{lem:weakconv2}
				Let $\mathfrak{P}_{N}$ be a sequence of point processes, and denote their $k-$ point correlation functions by $R_{k}^{N}$. Then $\mathfrak{P}_{N}$ converges weakly to a point process $\mathfrak{P}$, with $k-$ point correlation functions denoted by $R_{k}^{N}$ if:
				\begin{itemize}
					\item[1.] $R_{k}^{N} \to R_{k}$ pointwise almost everywhere as $N \to \infty$.
					
					\item[2.] For any compact set $\Omega \subset \mathbf{R}^{d}$ there holds
					\begin{equation}
						\sup_{N \in \mathbf{N}} \sum_{k=1}^{\infty} \frac{1}{k!} \int_{K} R_{k}^{N} \, \mathrm d y_{1}... \, \mathrm d y_{k} < \infty.
					\end{equation}
				\end{itemize}
			\end{lemma}
			
			Equipped with Lemma \ref{lem:weakconv2}, we will prove Theorem \ref{teo:poisson}, restated here for convenience. 
			
			\setcounter{theorem}{0}
			\begin{theorem}[Convergence to Poisson point process]
				Assume that $g$ is weakly interacting, and $V$ is admissible. Assume also that $\mathrm{ent}[\mu_{V}] < \infty$. Let $x^{*} \in \mathbf{R}^{d}$ and define the local point process $\Xi$ by 
				\begin{equation}
					\Xi := \sum_{i=1}^{N} \delta_{N^{\frac{1}{d}} \left( x_{i} - x^{*} \right)}.
				\end{equation}
				
				Then, for a.e. $x^{*}$, if $N^{-1} \ll \beta \ll N^{-\frac{1}{2}}$, $\Xi$ converges to a Poisson point process of intensity $\mu_{V}(x^{*})$.   
			\end{theorem}
			
			\begin{proof}
				We will assume that $N \beta$ is sufficiently large below.
				
				\textbf{Step 1:} Formula for the correlation functions.
				
				Let $k \leq \frac{1}{100} \sqrt{N}$. Define $Y_{k} := (y_{1}, ...y_{k}) \in \mathbf{R}^{d \times k}$, $x_{i} := x^{*} + N^{-\frac{1}{d}} (y_{i} - x^{*})$, and $X_{N} = (x_{1}, ... x_{N})$. Let $\tilde \emp_{N-k}$ be the empirical measure of $(x_{k+1},\ldots,x_N)$. Using Proposition \ref{prop:laplace} and Remark \ref{rem:corr}, the $k-$point correlation function of $\Xi$ can be written as  
				\begin{equation}
					\label{eq:corrfunc}
					\begin{split}
						R^{N}_{k}(Y_{k}) &= C(N,k)\frac{1}{Z_{N, \beta}} \int_{\R^{d(N-k)}} \exp \left( - \beta \mathcal{H}_{N}(X_{k}, x_{k+1},...x_{N}) \right) \, \mathrm d x_{k+1}... \, \mathrm d x_{N}\\
						&= C(N,k) \frac{1}{Z_{N, \beta}} \int_{\R^{d(N-k)}} \exp \left( - \beta \left( \sum_{1 \leq i \ne j \leq N} g(x_{i} - x_{j}) + N \sum_{i=1}^{N} V(x_{i}) \right) \right) \, \mathrm d x_{k+1}... \, \mathrm d x_{N}\\
						&= C(N,k)  \frac{Z_{N-k, \beta}}{Z_{N, \beta}} \exp \left( - \beta \left( \sum_{1 \leq i \ne j \leq k } g(x_{i} - x_{j}) + N \sum_{i=1}^{k} V(x_{i}) \right) \right) \\
						&\quad \times \mathbb{E}_{N-k, \beta} \left[ \exp\left( -(N-k) \beta \left( \sum_{i=1}^{k} 2h^{{\tilde {\rm emp}}_{N-k}} (x_{i}) \right) \right) \right]\\
						&= C(N,k)  \frac{Z_{N-k, \beta}}{Z_{N, \beta}} \exp \left( - \beta \left( \sum_{1 \leq i \ne j \leq k } g(x_{i} - x_{j}) + N \sum_{i=1}^{k} V(x_{i}) \right) \right) \\
						&\quad \times \left(  M_{N-k} \exp\left( -(N-k) \beta \left( \sum_{i=1}^{k} 2h^{\mu_{\theta}} (x_{i}) \right) \right) + A_{N-k} \right),
					\end{split}
				\end{equation}
				where $C(N,k) := \frac{N!}{(N-k)!} N^{-k}$, and $M_{N-k}$ and $A_{N-k}$ are as in Proposition \ref{prop:laplace} for $(x_{k+1},\ldots,x_N) \sim \P_{N-k,\beta}$.
				
				\textbf{Step 2:} Estimating the ratio of partition functions. 
				
				Integrating equation \eqref{eq:corrfunc}, and using that the k-th marginal is a probability measure, we get
				\begin{equation}
					\begin{split}
						\lefteqn{ \frac{Z_{N, \beta}}{Z_{N-k, \beta}}} \quad &\\
						&= C(N,k)  \int_{\mathbf{R}^{d \times k }} \exp \left( - \beta \left( \sum_{1 
							\leq i \ne j \leq k } g(x_{i} - x_{j}) + N \sum_{i=1}^{k} V(x_{i}) \right) \right) \\
						&\quad \times \left(  M_{N-k} \exp\left( -(N-k) \beta \left( 
						\sum_{i=1}^{k} 2h^{\mu_{\theta}} (x_{i}) \right) \right) + A_{N-k} \right) \, \mathrm d x_{1}... \, \mathrm d x_{k}\\
						&=  M_{N}'\int_{\mathbf{R}^{d \times k }} \exp \left( - N \beta \left( \sum_{i=1}^{k} V(x_{i}) + 2h^{\mu_{\theta}} (x_{i}) \right) \right) \, \mathrm d x_{1}... \, \mathrm d x_{k}+ A_{N}'
						\\
						&=  M_{N}' L_{\theta}^{k} + A_{N}', 
					\end{split}
				\end{equation}
				where
				\begin{equation}
					\label{eq:error'}
					\begin{split}
						|\log M_{N}'| &\leq C N^{\frac{1}{2}}\beta k \\
						|A_{N}'| &\leq C\exp \left( - \frac12 N \beta g(0) \right),
					\end{split}
				\end{equation}
				for some constant $C$ dependent only on $g$ and $V$. Here we used
				$$
				\left |  \sum_{1 \leq i < j \leq k} g(x_i - x_j) \right | \leq Ck^2 \leq C \sqrt{N} k,
				$$
				and that $\( \int_{\R^d} e^{-\beta V(x)} dx \)^k \leq C^k \leq \exp(\frac14 N \beta g(0))$.
				
				Recall that, by Remark \ref{rem:qual}, $L_{\theta}$ is of order $\exp \left( - 2 c_\theta N \beta \right)$. More specifically, 
				\begin{equation}
					\lim_{N \to \infty} \frac{\log L_{\theta}}{N \beta} =- 2 c_{\infty}. 
				\end{equation}
				By a Taylor expansion of $x \mapsto \frac{1}{x}$, we have that
				\begin{equation}
					\label{eq:ratio}
					\frac{Z_{N-k, \beta}}{Z_{N, \beta}} =  M'_{N} L_{\theta}^{-k} + A'_{N},
				\end{equation}
				for new terms $M'_N, A'_{N}$ satisfying equation \eqref{eq:error'} with a different constant $C$. 
				
				\textbf{Step 3}: Conclusion.
				
				\textbf{Subtep 3.1:} Pointwise convergence (item 1. of Lemma \ref{lem:weakconv2}).
				
				Let $k$ be fixed. Plugging in equation \eqref{eq:ratio} into equation \eqref{eq:corrfunc}, and using equation \eqref{foc} we obtain
				\begin{equation}
					\begin{split}
						R^{N}_{k}(Y_{k}) &= \left( M'_{N}  L_{\theta}^{-k} + A'_{N}\right) \exp \left( - \beta \left( \sum_{1 \leq i \ne j \leq k } g(x_{i} - x_{j}) + N \sum_{i=1}^{k} V(x_{i}) \right) \right)\\
						&\quad \times \left(  M_{N-k} \exp\left( -(N-k) \beta \left( \sum_{i=1}^{k} 2h^{\mu_{\theta}} (x_{i}) \right) \right) + A_{N-k} \right)\\
						&= M''_{N} L_{\theta}^{-k} \exp \left( -N \beta \left( \sum_{i=1}^{k} V(x_{i}) + 2h^{\mu_{\theta}} (x_{i}) \right)\right) +  A''_{N}\\
						&= M''_{N} \prod_{i=1}^{k} \mu_{\theta}(x_{i}) +  A''_{N}, 
					\end{split}
				\end{equation}
				where
				\begin{equation}
					\begin{split}
						|\log M''_{N}| &\leq CN^{\frac{1}{2}}\beta  \left(  k + \sum_{i=1}^{k}V(x_{i}) \right), \\
						|A''_{N}| &\leq C \exp \left( -  \frac12 N \beta   \left( g(0) + 2 \sum_{i=1}^k V(x_i) \right ) \right),
					\end{split}
				\end{equation}
				for some new constant $C > 0$.
				
				Using Lemma \ref{lem:convergence} (more specifically, equation \eqref{eq:microconv} with $\delta = N^{-\frac{1}{d}}$) we have that $\mu_{\theta}(x_{i})|_{K}$ converges strongly in $L^{1}$ (and modulo a subsequence, not relabelled, pointwise a.e.) to $\mu_{V}(x^{*})$. Using this, along with the hypothesis that $\beta N^{1/2} \ll 1$, we have that $R^{N}_{k}(Y_{k}) \to (\mu_{V}(x^{*}))^{k}$ for a.e. $Y_{k} \in \mathbf{R}^{d \times k}$. 
				
				\textbf{Substep 3.1:} Summability condition (item 2. of Lemma \ref{lem:weakconv2}). 
				
				Let $K \subset \mathbf{R}^{d}$ be a compact set. Then by Lemma \ref{lem:convergence}, there holds
				\begin{equation}
					\begin{split}
						\sup_{N \in \mathbf{N}} \sum_{k=1}^{\infty} \frac{1}{k!} \int_{K}R^{N}_{k}(Y_{k})  \, \mathrm d y_{1}... \, \mathrm d y_{k} &\leq  C\sum_{k=1}^{\infty} \frac{1}{k!} \left( |K|\mu_{V}(x^{*}) \right)^{k}
						< \infty.
					\end{split}
				\end{equation}
				
				Using Lemma \ref{lem:weakconv2}, we may conclude the proof.
				
				\section{Appendix: Point processes}
				
				In this appendix, we will review the theory on point processes necessary to prove Theorem \ref{teo:poisson}. This appendix is based on \cite{lambert2021poisson}. 
				
				\begin{definition}[Point process]
					A point process $\mathfrak{P}$ is a probability measure on the space of locally finite point configurations on $\mathbf{R}^{d}$. Alternatively, we may think of a point process $\mathfrak{P}$ as a probability measure on the space of Radon measures consisting of a sum of Dirac deltas. 
				\end{definition}
				
				\begin{definition}[Laplace functional]
					We associate to a point process $\mathfrak{P}$ its Laplace functional $\psi$, defined by
					\begin{equation}
						\psi(f) := \mathbb{E}_{\mathfrak{P}}\left[ \exp \left( - \int_{\mathbf{R}^{d}} f\,  \mathrm d X \right) \right],
					\end{equation}
					for a Borel-measurable $f: \mathbf{R}^{d} \to [0, \infty)$. 
				\end{definition}
				
				\begin{definition}[Correlation function]
					
					Given a point process, we define the $k-$point correlation functions $\{R_{k}\}_{k=1}^{\infty}$ by the condition that
					\begin{equation}
						\psi(f) := 1 + \sum_{i=1}^{\infty} \frac{1}{k!} \int_{\mathbf{R}^{d}} \prod_{i=1}^{k} \left( \exp(- f(x_{i})) -1 \right) R_{k}(x_{1},...x_{k}) \, \mathrm d x_{1}...\, \mathrm d x_{k}.
					\end{equation}
				\end{definition}
				
				\begin{remark}[Correlation function of a joint distribution]
					\label{rem:corr}
					If $\mathfrak{P}$ is supported on point configurations of $N$ points, and if these $N$ points have a symmetric joint distribution $\mathbb{P}$, then we can verify that
					\begin{equation}
						\psi(f) := 1 + \sum_{i=1}^{\infty}  {N \choose k} \int_{\mathbf{R}^{d}} \prod_{i=1}^{k} \left( \exp(- f(x_{i})) -1 \right)  \mathbb{P}(X_{N}) \, \mathrm d X_{N},
					\end{equation}
					which implies that the correlation functions are given by 
					\begin{equation}
						\label{eq:corr_sym}
						R_{k} (x_{1}, ...x_{k}) =  \frac{N!}{(N-k)!} \int_{\mathbf{R}^{d \times (N-k)}}  \mathbb{P}(x_{1},... x_{k}, x_{k+1}, ...x_{N}) \, \mathrm d x_{k+1} ... \, \mathrm d x_{N}
					\end{equation}
					for any $k \leq N$. 
					
				\end{remark}
				
				\begin{remark}[Correlation function of a Poisson point process]
					
					If $\mathfrak{P}$ is a Poisson process of intensity $\lambda$, then its Laplace functional is given by 
					\begin{equation}
						\psi(f) = \exp \left( \int_{\mathbf{R}^{d}} \left( \exp (-f(x)) -1  \right) \lambda(x) \, \mathrm d x \right). 
					\end{equation}
					
					Furthermore, the correlation functions of $\mathfrak{P}$ are given by
					\begin{equation}
						R_{k}(x_{1},...x_{k}) = \prod_{i=1}^{k} \lambda(x_{i}). 
					\end{equation}
				\end{remark}
				
				\begin{definition}[Weak convergence]
					
					Let $\mathfrak{P}_{N}$ be a sequence of point processes, and denote their Laplace functional by $\psi_{N}$. Then $\mathfrak{P}_{N}$ converges weakly to a point process $\mathfrak{P}$ with Laplace functional $\psi$ if for any continuous and compactly supported $f: \mathbf{R}^{d} \to \mathbf{R}^{+}$ there holds
					\begin{equation}
						\lim_{N \to \infty} \psi_{N}(f) = \psi(f).
					\end{equation}
				\end{definition}
				
				\begin{lemma}[Weak convergence and correlation functions]
					\label{lem:weakconv}
					Let $\mathfrak{P}_{N}$ be a sequence of point processes, and denote their $k-$ point correlation functions by $R_{k}^{N}$. Then $\mathfrak{P}_{N}$ converges weakly to a point process $\mathfrak{P}$, with $k-$ point correlation functions denoted by $R_{k}^{N}$ if:
					\begin{itemize}
						\item[1.] $R_{k}^{N} \to R_{k}$ pointwise almost everywhere as $N \to \infty$.
						
						\item[2.] For any compact set $\Omega \subset \mathbf{R}^{d}$ there holds
						\begin{equation}
							\sup_{N \in \mathbf{N}} \sum_{k=1}^{\infty} \frac{1}{k!} \int_{K} R_{k}^{N} \, \mathrm d y_{1}... \, \mathrm d y_{k} < \infty.
						\end{equation}
					\end{itemize}
				\end{lemma}

			\end{proof}
			\section{Acknowledgements}
			
			DPG acknowledges the support of the ERC Starting Grant ``Bridging Scales in Random Materials" ERC StG RandSCALES 948819. ET was supported by NSF grant DMS-2303318.

			\bibliographystyle{plain}
			\bibliography{bibliography.bib}
			
			\vskip .5cm
			\noindent
			\textsc{David Padilla-Garza}\\
			Institute of Science and Technology Austria. \\
			Email: {David.Padilla-Garza@ist.ac.at}.
			\vspace{.2cm}
			
			\vskip .5cm
			\noindent
			\textsc{Luke Peilen}\\
			Department of Mathematics, Temple University. \\
			Email: {luke.peilen@temple.edu}.
			\vspace{.2cm}

			\vskip .5cm
			\noindent
			\textsc{Eric Thoma}\\
			Department of Mathematics, Stanford University. \\
			Email: {thoma@stanford.edu}.
			\vspace{.2cm}
			
		\end{document}